\documentclass[12pt,twoside]{article}

\usepackage{amssymb,amsthm,amsmath,amsfonts,bbold}

\def\argmax{\mathop{\rm argmax}} 
\def\rbox{\begin{flushright} $\square $ \end{flushright} }

\begin{document}

\emph{Theoretical Mathematics \& Applications, vol. x, no. x, 2013, x-x}

\emph{ISSN: 1792-9687 (print), 1792-9709 (online)}

{Scienpress Ltd, 2013}

\centerline{}

\centerline{}

\centerline {\Large{\bf A Consistency Result for Bayes Classifiers}}

\centerline{}

\centerline{\Large{\bf  with Censored Response Data}}

\centerline{}
\centerline{\bf {Priyantha Wijayatunga\footnote{Department of Statistics,  Ume\r{a} School of Business and Economics, Ume\r{a}  University, Ume\r{a} 90187, Sweden. e-mail: priyantha.wijayatunga@stat.umu.se} and
Xavier de Luna\footnote{Department of Statistics,  Ume\r{a} School of Business and Economics, Ume\r{a}  University, Ume\r{a} 90187, Sweden. e-mail: xavier.deluna@stat.umu.se
~~~~~~ \hspace{2.9cm }
\newline Article Info: {\it Received} : January 10, 2011.~~{\it Published online} : July 17, 2011}}}

\centerline{}

\newtheorem{Theorem}{\quad Theorem}[section]

\newtheorem{Assumption}[Theorem]{\quad Assumption}

\newtheorem{Definition}[Theorem]{\quad Definition}

\newtheorem{Corollary}[Theorem]{\quad Corollary}

\newtheorem{Lemma}[Theorem]{\quad Lemma}

\newtheorem{Example}[Theorem]{\quad Example}

\def\iDelta{\mathit{\Delta}}
\def\iTheta{\mathit{\Theta}}
\def\bix{\boldsymbol{x}}
\def\biX{\boldsymbol{X}}
\def\biz{\boldsymbol{z}}
\def\biE{\boldsymbol{E}}
\def\biP{\boldsymbol{P}}
\def\biD{\boldsymbol{D}}
\def\bid{\boldsymbol{d}}
\def\bione{\boldsymbol{1}}
\def\bizero{\boldsymbol{0}}
\def\bitheta{\boldsymbol{\theta}}
\def\bialpha{\boldsymbol{\alpha}}

\def\bbN{{\mathbb N}}
\def\bbR{{\mathbb R}}
\def\bbone{\text{\textbb 1}}

\begin{abstract}
Naive Bayes classifiers have proven to be useful in many prediction problems with complete training data. Here we consider the situation where a naive Bayes classifier is trained with data where the response is right censored. Such prediction problems are for instance encountered in profiling systems used at National Employment Agencies. In this paper we propose the maximum collective conditional likelihood estimator for the prediction and show that it is strongly consistent under the usual identifiability condition. 
\end{abstract}
\vspace{0.76cm}
{\bf Mathematics Subject Classification :} xxxxx 

\section{Introduction} Naive Bayes classifiers have proven to be useful in many prediction problems with complete training data. Here we consider the situation where a naive Bayes classifier is trained with data where the response is right censored. Such prediction problems are for instance encountered in profiling systems used at National Employment Agencies.  A profiling system provides predictions of unemployment duration based on individual characteristics. To train such a system register data on individuals is used where unemployment durations as well as demographic and socio-economics information is recorded. Unemployment duration is then typically censored by the end of the observation period as well as exit from unemployment due to other reasons than employment, typically entrance into educational programs \cite{NP08}. Naive Bayes classifiers have proven useful in many prediction problems with complete data \cite{WMN06}, \cite{WBW05} and references therein. In this paper for this censored response case we propose the maximum collective conditional likelihood estimator and show that it is strongly consistent under the usual identifiability condition \cite{WM07} whose notation is used in our proof. 

Formally, consider a class variable $X_0$ $-$say unemployment duration$-$ and a $n-$attribute random variable
vector $X_{[n]}=(X_1, \ldots, X_n)$ $-$individual features, where all the variables are discrete
and finite. Note that contitunous variables may be discretized making this framework more general. We assume thus that the state space of $X_i$ is
$\mathcal{X}_i = \{1, \ldots, r_i\}$. We further assume that
$(X_0,X_{[n]})$ forms a naive Bayesian network so that their
joint density and conditional density of $X_{0}$ given $X_{[n]}=x_{[n]}$ are as follows:
\begin{align*}
p(x_0,x_1,\ldots,x_n) &= p(x_0) \prod_{i=1}^n p(x_i \mid x_0), \\
p(x_0\mid x_{[n]})    &= \frac{p(x_0) \prod_{i=1}^n p(x_i\mid x_0)}{\sum_{x'_0} p(x'_0) \prod_{i=1}^n p(x_i\mid x'_0)}.
\end{align*}

Let the parameter space be $\iTheta = \iDelta_{r_0}\times \iDelta_{r_1}^{r_0}\times \cdots \times \iDelta_{r_n}^{r_0}$
where $\iDelta_t = \{ (p_1,\ldots,p_{t-1}): 0\leq p_i, \sum_{i=1}^{t-1} p_i \leq 1 \}$ and for $\iDelta_a^b$, where $b$ is the usual power.
The interior of $\iTheta$ is denoted by $\iTheta^o$. In the following, we always
assume that the true parameter is an element of $\iTheta^o$. Note that,
if this is not the case then the naive  Bayesian network is degenerated
in the sense that some variables (if binary) may vanishes as their state space
will be reduced to singletons or some may appear with reduced state
spaces.   Then, the parameter of the naive Bayesian model is $\theta =
(\theta_0, \theta_1, \ldots, \theta_n)\in
\iTheta$, where  $\theta_0=(\theta_{x_0=1}, \ldots,
\theta_{x_0=r_0-1})$ and $\theta_i = (\theta_{i \mid x_0=1}, \ldots,
\theta_{i \mid x_0=r_0})$ such that $\theta_{i \mid
x_0}=(\theta_{x_i=1 \mid x_0},\dots, \theta_{x_i=r_i-1 \mid x_0})$  for $i=1,
\dots, n$.  Since we are working on the non-Bayesian case we have

\begin{align*}
p_{\theta}(x_0)         &= \theta_{x_0}               \qquad \text{ if $x_0=1,\ldots,r_0-1$}, \\
p_{\theta}(r_0)         &= 1-\sum_{x_0=1}^{r_0-1} \theta_{x_0}, \\
p_{\theta}(x_i\mid x_0) &= \theta_{x_i\mid x_0}       \qquad \text{ if $x_i=1,\ldots,r_i-1$}, \\
p_{\theta}(r_i\mid x_0) &= 1-\sum_{x_i=1}^{r_i-1} \theta_{x_i\mid x_0}.
\end{align*}
Note that the above marginal and conditional densities are
over-parameterized, i.e., when we write, for example, the density of
$X_0$, $p_\theta(x_0)$ the parameter $\theta$ also contains irrelevant
components in addition to relevant ones
to determine the probabilities of $X_0$.

Suppose we have a random sample with $N=N_1+N_2$ number of data cases on the random variable vector $(X_0,X_1,\ldots,X_n)$, and denote 
$$D=
\big\{(x_0^{(1)},x_1^{(1)},\ldots,x_n^{(1)}),\ldots,(x_0^{(N)},x_1^{(N)},\ldots,
x_n^{(N)}) \big\}.$$ 
Let the first $N_1$ cases to be fully, and for the remaining $N_2$ cases to be right censored in $X_0$. In the example of unmeployment duration, where $X_0$ denotes the duration of an unemployment spell for an individual and $(X_1,\ldots,X_n)$ a suitable vector of features/covariates, (right) censoring of $X_0$ may be due to the end of study or entrance in educational programs for unemployed. In the sequel, by ``censored'' it is meant right censored response.  

The compound collective conditional likelihood (CCCL) of $\theta$ given the data $D$ is defined as
\begin{eqnarray*}
CCCL_N(\theta)&=& \prod_{j=1}^{N_1} p_{\theta}(x_0^{(j)} \mid x_{[n]}^{(j)}) \prod_{k=N_1+1}^{N}P_{\theta}(X_0 > x_0^{(k)} \mid x_{[n]}^{(k)}) \\
             &=& \prod_{j=1}^{N_1} p_{\theta}(x_0^{(j)} \mid x_{[n]}^{(j)}) \prod_{k=N_1+1}^{N}\sum_{x_0=x_0^{(k)}+1}^{r_0}p_{\theta}(x_0 \mid x_{[n]}^{(k)}) \\
             &=& \prod_{j=1}^{N_1} \frac{p_{\theta}(x_0^{(j)}) \prod_{i=1}^n p_{\theta}(x_i^{(j)} \mid x_0^{(j)})}
                  {\sum_{x'_0} p_{\theta}(x'_0) \prod_{i=1}^n p_{\theta}(x_i^{(j)} \mid x'_0)}   \\
             & & \times \prod_{k=N_1+1}^{N} \sum_{x_0=x_0^{(k)}+1}^{r_0} \frac{p_{\theta}(x_0) \prod_{i=1}^n p_{\theta}(x_i^{(k)} \mid x_0)}
                  {\sum_{x'_0} p_{\theta}(x'_0) \prod_{i=1}^n p_{\theta}(x_i^{(k)} \mid x'_0)} \\
             &= & \prod_{j=1}^{N_1} \frac{p_{\theta}(x_0^{(j)}) \prod_{i=1}^n p_{\theta}(x_i^{(j)} \mid x_0^{(j)})}
                  {\sum_{x'_0} p_{\theta}(x'_0) \prod_{i=1}^n p_{\theta}(x_i^{(j)} \mid x'_0)}   \\
   & & \times \prod_{k=N_1+1}^N \Bigg\{ \frac{P_{\theta}(X_0=x_0^{(k)}+1) \prod_{i=1}^n p_{\theta}(x_i^{(k)} \mid X_0=x_0^{(k)}+1)}
                  {\sum_{x'_0} p_{\theta}(x'_0) \prod_{i=1}^n p_{\theta}(x_i^{(k)} \mid x'_0)}  \\
   & & + ... + \frac{P_{\theta}(X_0=r_0) \prod_{i=1}^n p_{\theta}(x_i^{(k)} \mid X_0=r_0)}
                  {\sum_{x'_0} p_{\theta}(x'_0) \prod_{i=1}^n p_{\theta}(x_i^{(k)} \mid x'_0)} \Bigg\} \\
   &=& CCL_{N}^1(\theta)+... +CCL_{N}^M(\theta), 
\end{eqnarray*}
where $M=\prod_{j=N_1+1}^N (r_0-x_0^{(j)})$. Thus, the CCCL is a sum of $M$ collective conditional likelihoods.

The maximum compound collective conditional likelihood estimator (MCCCLE) is then
\begin{gather*}
\hat{\theta}_N = \argmax_{\theta \in \iTheta} \Big\{ CCL_{N}^1(\theta)+ ...+ CCL_{N}^M(\theta) \Big\}.
\end{gather*}
MCCCLE has no closed form expression in general and need to be solved numerically. We will need below the following MCCLEs
\begin{eqnarray*}
  \hat{\theta}_N^i = \argmax_{\theta \in \iTheta} CCL_{N}^i(\theta) \ \ \ \textrm{ for } i=1,...,M.
\end{eqnarray*}

\section{Strong consistency of MCCCLE}

In this section, we give a proof for the strong consistency of MCCCLE. First, we need the following identifiability assumption as usual in maximum likelihood theory. 

\begin{Assumption} \label{identifiability} \emph{(Identifiability Condition)}
If $p_{\theta}(x_0\mid x_{[n]})=p_{\theta'}(x_0\mid x_{[n]})$
for all $x_0$ and $x_{[n]}$ then $\theta=\theta'$.
\end{Assumption}
This condition requires that $\theta$ should be uniquely
determined by the corresponding density $p_\theta(.\mid .)$. As shown earlier \cite{WM07} for MCCLE, this does not always hold.

The collective conditional likelihood function for the naive Bayes network model with complete data is a concave down function in the parameters \cite{RGMT05}. As noticed above, the compound collective conditional likelihood function with censored response is a sum of collective likelihood function.

\begin{Lemma}
The compound collective conditional likelihood function $CCCL_N(\theta)$ defined above is concave down in $\theta$.
\end{Lemma}

\noindent
\emph{\bf Proof:} First note that for any $\theta \in \iTheta$, $p_\theta(x_0 \mid x_{[n]})>0$, and, therefore, so do the collective conditional likelihood functions composing $CCCL_N$. Furthermore, the sum of two convex functions with same support is convex and so do the sum of any number of convex functions with same domain, thereby yielding the result.
\rbox
\begin{Lemma}
If $f$ and $g$ are two convex functions on the same domain with their global minima at $x_1$ and $x_2$ respectively, then $f+g$ has its global minima at $tx_1+(1-t)x_2$ for some $t \in [0,1]$.
\end{Lemma}
\noindent
\emph{\bf Proof}: 
If $x_1=x_2$ then the result holds for $t=0$ or $t=1$. Consider the case where $x_1 \neq x_2$. Then $\dot{f}(x_1)+\dot{g}(x_1) < 0 $ and  $\dot{f}(x_2)+\dot{g}(x_2) > 0 $ or vice versa. Since $f+g$ is a convex function (sum of two convex functions with same domain), there  must be a point $x_3$ such that $\dot{f}(x_3)+\dot{g}(x_3) = 0$, where $x_3=tx_1+(1-t)x_2$ for some  $t \in [0,1] $.  
\rbox

If  $M=2$ then, since  $CCL_N^1(\theta)$ and $CCL_N^2(\theta)$ are concave down functions having their  maxima at $\hat{\theta}_N^1$ and $\hat{\theta}_N^2$ respectively, $CCL_N^1(\theta)+ CCL_N^2(\theta)$ which is also concave down has its maximum at  $\hat{\theta}_N:=t_D \hat{\theta}_N^1 + (1-t_D) \hat{\theta}_N^2 $ where $t_D $ is a vector of the same length as $\theta^*$, whose components are in $ [0,1] $, and which is dependent on the data $D$. Note that both $\hat{\theta}_N^1$ and $\hat{\theta}_N^2$ are consistent estimates for $\theta^*$ under Assumption \ref{identifiability} (since they are MCCLEs, see \cite{WM07}).

\begin{eqnarray}
  P_{\theta^*} \Big\{  \lim_{N_1 \rightarrow \infty } \hat{\theta}_N^1  &=& \theta^*  \Big\} =1 \label{theta1.eq}\\
 P_{\theta^*} \Big\{  \lim_{N_1 \rightarrow \infty } \hat{\theta}_N^2  &=& \theta^*  \Big\} =1 \label{theta2.eq}
\end{eqnarray}
In the sequel, we write $\theta=0$ for $\theta \in \iTheta^0$ to mean that  all the components of $\theta$ are zeros and similarly for any inequality on $\theta$. Now we can prove the strong consistency of MCCCLE.

\begin{Theorem} \label{strong.consistency}
Under Assumption \ref{identifiability}, MCCCLE $\hat{\theta}_N$ is strongly consistent as as follows,
\begin{eqnarray}
 P_{\theta^*} \Big\{  \lim_{N_1 \rightarrow \infty } \hat{\theta}_{M+N_1}  &=& \theta^*  \Big\} =1, \forall M.
\end{eqnarray}
\end{Theorem}

\noindent
\emph{\bf Proof:} We prove the result by induction. Let $M=2$, 
then $\hat{\theta}_N = t_D \hat{\theta}_N^1+ (1-t_D) \hat{\theta}_N^2 = \hat{\theta}_N^2 + t_D ( \hat{\theta}_N^1 -  \hat{\theta}_N^2 )$. 
By (\ref{theta1.eq}) and (\ref{theta2.eq}) we have
\begin{eqnarray}
 P_{\theta^*} \Big\{ \lim_{N_1 \rightarrow \infty } \hat{\theta}_N^1  = \theta^* \cap \lim_{N_1 \rightarrow \infty } \hat{\theta}_N^2  = \theta^*    \Big\} =1. \\
\end{eqnarray}
Since  $ \hat{\theta}_N^1,  \hat{\theta}_N^2$ and  $\theta^*$ are finite and $0 \leq t_D \leq 1$ (therefore $0 \leq \limsup_{N \rightarrow \infty}  t_D \leq 1$) we can write
\begin{eqnarray}
 P_{\theta^*} \Big\{ \limsup_{N_1 \rightarrow \infty } t_D(\hat{\theta}_N^1  - \hat{\theta}_N^2)  = 0    \Big\} &=&1 \\
P_{\theta^*} \Big\{ \limsup_{N_1 \rightarrow \infty } \hat{\theta}_N^2 + t_D(\hat{\theta}_N^1  - \hat{\theta}_N^2)  = \theta^*    \Big\} &=& 1 \\
P_{\theta^*} \Big\{ \limsup_{N_1 \rightarrow \infty } \hat{\theta}_N  = \theta^*   \Big\} &=& 1
\end{eqnarray}
Similarly we can write 
\begin{eqnarray}
P_{\theta^*} \Big\{ \liminf_{N_1 \rightarrow \infty } \hat{\theta}_N  = \theta^*   \Big\} &=& 1.
\end{eqnarray}
Hence
\begin{eqnarray}
P_{\theta^*} \Big\{ \lim_{N_1 \rightarrow \infty } \hat{\theta}_N  = \theta^*   \Big\} &=& 1
\end{eqnarray}
Now assume that for $M>2$,  $\hat\theta_N := w_D^1 \hat{\theta}_N^1+ ...+w_D^M \hat{\theta}_N^M$ where $w_D^1 +...+w_D^M=1$, the maximizer of $CCCL_N(\theta)$, is a consistent estimator of $\theta^*$. Then, for the case of $M+1$, assume for simplicity that the additional new censored observation is  $x_0^{(N+1)}=r_0-2$. Then,
\begin{eqnarray*}
CCCL_{N+1}(\theta)&=&  \Big\{ CCL_{N}^1(\theta)+... +CCL_{N}^M(\theta)  \Big\} \\
   & & \times  \Bigg\{ \frac{P_{\theta}(X_0=r_0-1) \prod_{i=1}^n p_{\theta}(x_i^{(N+1)} \mid X_0=r_0-1)}
                  {\sum_{x'_0} p_{\theta}(x'_0) \prod_{i=1}^n p_{\theta}(x_i^{(N+1)} \mid x'_0)}  \\
   & & + \frac{P_{\theta}(X_0=r_0) \prod_{i=1}^n p_{\theta}(x_i^{(N+1)} \mid X_0=r_0)}
                  {\sum_{x'_0} p_{\theta}(x'_0) \prod_{i=1}^n p_{\theta}(x_i^{(N+1)} \mid x'_0)} \Bigg\}. 
\end{eqnarray*}
Let us rewrite this (with obvious new notation) 
\begin{eqnarray*}
CCCL_{N+1}(\theta)&=& \Big\{ CCL_{N+1}^{a,1}(\theta)+... +CCL_{N+1}^{a,M}(\theta)  \Big\} \\
               & & + \Big\{ CCL_{N+1}^{b,1}(\theta)+... +CCL_{N+1}^{b,M}(\theta)  \Big\}
\end{eqnarray*}
Now denote $\hat{\theta}_{N+1}^a := w_D^{a,1} \hat{\theta}_{N+1}^{a,1}+ ...+w_D^{a,M} \hat{\theta}_{N+1}^{a,M}$ and  $\hat{\theta}_{N+1}^b := w_D^{b,1} \hat{\theta}_{N+1}^{b,1}+ ...+w_D^{b,M} \hat{\theta}_{N+1}^{b,M}$ 
where $\sum_j^M w_D^{i,j}=1$ for $i=a,b$, the maximizers of the first and second sums of CCLs respectively. By assumption they are consistent estmators of $\theta^*$. Now similarly to the case $M=2$, we can write $\hat{\theta}_{N+1} := u_D\hat{\theta}_{N+1}^a+(1-u_D)\hat{\theta}_{N+1}^b$, where $0 \leq u_D \leq 1$, the maximizer of $CCCL_{N+1}(\theta)$, and show that it is strongly consistent for $\theta^*$. 
\rbox

\begin{Corollary} \label{cor.strong.consistency}
$p_{\hat{\theta}_N}(x_0\mid x_{[n]})$ is strongly consistent estimator of $p_{\theta^*}(x_0\mid x_{[n]})$
for each $x_{[n]}$.
\end{Corollary}

\noindent
\emph{\bf Proof:}
Immediate from the theorem since the densities $p_{\theta}(x_0\mid
x_{[n]})$ for all $x_{[n]}$ are rational functions of the parameter
which have no poles in $\iTheta^o$.
\rbox

{\bf ACKNOWLEDGEMENTS.} The authors gratefully acknowledge the financial support of the Swedish Research Council (through the Swedish Initiative for Microdata Research in the Medical and Social Sciences (SIMSAM) and Ageing and Living Condition Program)  and Swedish Research Council for Health, Working Life and Welfare (FORTE).

\end{document}